\documentclass[10pt,leqno]{article}
\usepackage{graphics}
\usepackage{amsmath,amssymb}

\newcommand{\om}{\omega}

\newcommand{\la}{\lambda}
\newcommand{\bt}{\beta}

\newcommand\R{\mathbb{R}}
\newcommand\C{\mathbb{C}}
\newcommand\huno{H^{1}}

\newcommand{\rn}{\R^{n}}

\newcommand{\lapl}{\Delta}

\newtheorem{theorem}{Theorem}
\newtheorem{remark}[theorem]{Remark}
\newtheorem{lemma}[theorem]{Lemma}
\newtheorem{proposition}[theorem]{Proposition}
\newtheorem{corollary}[theorem]{Corollary}

\newcommand{\bos}{\begin{remark}\rm}
\newcommand{\eos}{\end{remark}}
\newcommand{\bte}{\begin{theorem}}
\newcommand{\ete}{\end{theorem}}
\newcommand{\bpr}{\begin{proposition}}
\newcommand{\epr}{\end{proposition}}
\newcommand{\bco}{\begin{corollary}}
\newcommand{\eco}{\end{corollary}}
\newcommand{\ble}{\begin{lemma}}
\newcommand{\ele}{\end{lemma}}
\newcommand{\bdm}{\begin{displaymath}}
\newcommand{\edm}{\end{displaymath}}
\newcommand{\beq}{\begin{equation}}
\newcommand{\eeq}{\end{equation}}

\newcommand{\dimo}{{\noindent{\bf Proof.}\quad}}
\newcommand{\finedim}{{\unskip\nobreak\hfil\penalty50
 \hskip2em\hbox{}\nobreak\hfil\mbox{\rule{1ex}{1ex} \qquad}
  \parfillskip=0pt \finalhyphendemerits=0\par\medskip}}

\begin{document}

\title{On the blow-up threshold for weakly coupled nonlinear Schr\"odinger equations}
\author{Luca Fanelli\thanks{Partially supported by MIUR project
{\it Buona positura e stime di decadimento per equazioni dispersive e sistemi iperbolici}.}  \ \ 
 \& Eugenio Montefusco \thanks{Partially supported by MIUR project
{\it Metodi Variazionali ed Equazioni Differenziali non lineari}.}}
\date{\small Dipartimento di Matematica,\\
{\it Sa\-pien\-za} Universit\`a di Roma,\\
Piazzale Aldo Moro 5, 00185 Roma, Italy.}
\maketitle

\begin{abstract}
We study the Cauchy problem for a system of two coupled nonlinear
focusing Schr\"odinger equations arising in nonlinear optics.
We discuss when the solutions are global in time or blow-up in finite
time. Some results, in dependence of the data of the problem, are proved;
in particular we give a bound, depending on the coupling parameter, for the
blow-up threshold.

\end{abstract}

\section{Introduction}\label{sec.introd}

In this paper we consider the following Cauchy problem for two
coupled nonlinear Schr\"odinger equations
\begin{equation}\label{eq.schr1}
  \begin{cases}
  i\phi_t+\lapl\phi+\left(|\phi|^{2p}+\beta|\psi|^{p+1}|\phi|^{p-1}\right)\phi=0 & \\
  i\psi_t+\lapl\psi+\left(|\psi|^{2p}+\beta|\phi|^{p+1}|\psi|^{p-1}\right)\psi=0 & \\
  \phi(0,x)=\phi_0(x)\qquad\psi(0,x)=\psi_0(x),
  \end{cases}
\end{equation}
where $\phi,\psi:\R\times\R^n\to\C$, $\phi_0,\psi_0:\R^n\to\C$,
$p\geq0$ and $\beta$ is a real positive constant.

This kind of problems arises as a model for propagation of polarized
laser beams in birefringent Kerr medium in nonlinear optics (see,
for example, \cite{men, b, fp, lkr} and the references therein for a
complete discussion of the physics of the problem). The two
functions $\phi$ and $\psi$ are the components of the slowly varying
envelope of the electrical field, $t$ is the distance in the
direction of propagation, $x$ are the orthogonal variables and
$\lapl$ is the diffraction operator. The case $n=1$ corresponds to
propagation in a planar geometry, $n=2$ is the propagation in a bulk
medium and $n=3$ is the propagation of pulses in a bulk medium with
time dispersion (in this case $x$ includes also the time variable).

The focusing nonlinear terms in \eqref{eq.schr1} describe the dependence 
of the  refraction index of the material on the electric field intensity 
and the birefringence effects. The parameter $\bt>0$ has to be interpreted 
as the {\it birefringence} intensity and describes the coupling between the 
two components of the electric field envelope.
The case $p=1$ (i.e. cubic nonlinearities in \eqref{eq.schr1}) is 
known as Kerr nonlinearity in the physical literature.

We are interested in a slightly more general model, in order to
cover the physical cases and to discuss some results about Cauchy
problem \eqref{eq.schr1} from a more general point of view.

The aim of this paper is to study the $H^1\times H^1$ well-posedness 
of problem \eqref{eq.schr1}, with respect to the nonlinearity, in analogy 
with the case of the single focusing nonlinear Schr\"odinger equation
\begin{equation}\label{eq.nls}
  \begin{cases}
  i\psi_t+\Delta \psi+|\psi|^{2p}\psi=0
  \\
  \psi(0,x)=f(x),
  \end{cases}
\end{equation}
for $\psi:\R^{1+n}\to\C$ and $f:\R^n\to\C$. 

The well known results for the single equation can be summarized as follows. 
By standard scaling arguments it possible to claim that the critical exponent 
for the $H^1$ local well-posedness of \eqref{eq.nls} is $p=2/(n-2)$ 
(see \cite{caz}). Indeed, contraction techniques based on Strichartz 
estimates (see \cite{gv2}, \cite{kt}), permit to prove that \eqref{eq.nls} is 
locally well-posed in $H^1$ for $p<2/(n-2)$ (see \cite{gv1}, \cite{caz}). 

To pass from local to global well-posedness, it is natural to introduce
the energy function given by
\begin{equation*}
   E(t)=
   \frac12\|\nabla\psi\|_2^2-
   \frac1{2p+2}\|\psi\|_{2p+2}^{2p+2},
\end{equation*}
that is conserved along any solution $\psi$ of \eqref{eq.nls}. 
For $p<2/n$ the unique local $H^1$ solution can be extended globally 
in time by a continuation argument.
In the {\it critical case} $p=2/n$, we can also extend local solutions 
to global ones, provided the initial data are not too large in the $L^2$. 
Finally, for $2/n\leq p<2/(n-2)$ without restriction on the data, it is 
possible to prove that the $L^2$ norm of the gradient, in general, blows up 
in a finite time (see e.g. the original work \cite{gl} or \cite{caz}).  

By a physical point of view it is very interesting to determine the threshold for 
the initial mass of the wave packet, that is the $L^2$ norm of the initial datum, 
which separates global existence and blow-up in the critical case. 
We recall that $\psi=e^{it}u(x)\in H^1$ is a {\it ground state solution} for 
\eqref{eq.nls} if  $u$ is a nonzero critical point of the action functional
\bdm
  A(u)= E(u)+\frac{1}{2}\|u\|_2^2=\frac12\left(\|\nabla u\|_2^2
  +\|u\|_2^2\right) -\frac1{2+4/n}\|u\|_{2+4/n}^{2+4/n},
\edm 
having the smallest action level; moreover it follows that $u$ is a positive and radial solution of
\beq\label{eq.ellsing}
-\Delta u+ u=|u|^{4/n}u.
\eeq
In \cite{w}, Weinstein proved that if the initial mass is smaller than a constant
$C_n$, depending only on the space dimension $n$, than there exists a unique 
global $H^1$ solution; moreover $C_n$ is the $L^2$ norm of any ground state 
solutions of \eqref{eq.nls} and can be numerically estimated. Moreover we 
want to point out that this kind of phenomena for the single equation present other 
kind of universality properties related, for example, to the blow-up profile
(see \cite{b, mr, mr2} and the references therein). 

Our main goal is to state the analogous result for the coupled system 
\eqref{eq.schr1}. The critical exponent for the local $H^1\times H^1$ 
well-posedness has to be again $p=2/(n-2)$; so for $p<2/(n-2)$ it is possible to
prove that \eqref{eq.schr1} possesses a unique local solution (see Remark 
$4.2.13$ in \cite{caz} and Section \ref{sec.glob} below). 
The natural energy for \eqref{eq.schr1} is the following:
\begin{equation}\label{eq.energy}
   E(t)= \frac12\left(\|\nabla\phi\|_2^2+
   \|\nabla\psi\|_2^2\right)-\frac1{2p+2}\left(\|\phi\|_{2p+2}^{2p+2}
   +2\bt\|\phi\psi\|_{p+1}^{p+1}+ \|\psi\|_{2p+2}^{2p+2}\right).
\end{equation}
Also here it is possible to prove that $E(t)$ is conserved (see Section 
\ref{sec.glob}); hence the same techniques for the single equation can be 
applied to extend local solution to global ones. Now we can state our first
result.

\begin{theorem}\label{thm.existence}
Assume that $p<2/n$. Then the Cauchy problem \eqref{eq.schr1} is globally 
well posed in $H^1\times H^1$, i.e. for any $(\phi_0,\psi_0)\in H^1\times H^1$
there exists a unique solution $(\phi,\psi)\in C\left(\R;H^1\times H^1\right)$.
\end{theorem}

Also for system \eqref{eq.schr1}, the ground state solutions play a crucial 
role in the dynamics. In this case they are solutions of the form 
$(\phi,\psi)=e^{it}\left(u(x),v(x)\right)$, where the functions 
$u$ and $v$ have to be a least action solution of a elliptic system
(see \eqref{eq.ellittico} below). 

Since the birefringence tends to split a pulse into two pulses in two 
different polarization directions, the properties of the ground state solutions 
of \eqref{eq.schr1} depend strongly on the coupling parameter. 
If $\beta$ is sufficiently small, that is the interaction is weak, 
any ground state is a {\it scalar} solution, i.e. one of the two components 
is zero. On the other hand when the birefringence is strong, $\bt\gg1$,
we have {\it vector} ground states, i.e. all the components are 
distinct from zero (see \cite{ac, mmp}). 
This suggest that the also the blow-up phenomena, in the critical case, 
should depend on the parameter $\bt$. It is natural to claim that in a weak
interaction regime the behaviour has to be exactly the same of the single equation.
Otherwise if $\beta\gg1$, we expect that the analogous of the Weinstein threshold
$C_n$ should depend on $\bt$ also.
These claims are proved in the following main theorem. 

\begin{theorem}\label{thm.main}
Assume that $p=2/n$. Then there exists a constant $C=C_{n,\beta}$
such that the Cauchy problem \eqref{eq.schr1} is globally well posed
in $H^1\times H^1$ if
\begin{equation*}
\|\phi_0\|^2_2+\|\psi_0\|^2_2< C.
\end{equation*}
Moreover there exists a pair $(\phi_0,\psi_0)$ such that 
$\|\phi_0\|_2^2+\|\psi_0\|_2^2=C_{n,\bt}$ and the corresponding
solution blows up in a finite time.
The constant $C_{n,\bt}$ has the following behaviour
\beq\label{eq.cci}
\begin{cases}
C_{n,\bt}=C_{n} & \hbox{if } \bt\leq 2^{2/n}-1, \\
C_{n,\bt}\geq C_{n} \dfrac{(1+\bt)}{2^{2/n}} & \hbox{if } \bt\geq 2^{2/n}-1,
\end{cases} 
\eeq
where $C_{n}$ is the blow-up threshold of a single equation.
\ete

\bos In the supercritical case the solution of the Cauchy problem for 
\eqref{eq.schr1} exists locally in time, by the results in \cite{caz}. 
It is possible to prove that the solution exists globally in time if the 
assumption $\|\phi_0\|^2,\|\psi_0\|^2\ll 1$ is satisfied (see Theorem 
$6.1.1$ in \cite{caz}). 
\eos

\bos\label{rm.cubico}
As observed above, the Kerr nonlinearities (corresponding to $p=1$)
are physically relevant; in this case the system \eqref{eq.schr1} becomes
\begin{equation}\label{eq.schr2}
  \begin{cases}
  i\phi_t+\lapl\phi+ \left(|\phi|^{2}+\beta|\psi|^{2}\right)\phi=0 & \\
  i\psi_t+\lapl\psi+\left(|\psi|^{2}+\beta|\phi|^{2}\right)\psi=0.
  \end{cases}
\end{equation}

The above results (Theorems \ref{thm.existence} and \ref{thm.main}) can be summarized
in the following way:
\begin{itemize}
\item[i)] if $n=1$ the Cauchy problem \eqref{eq.schr2} is globally
  (in time) well posed in $H^1\times H^1$,
\item[ii)] if $n=2$ the cubic nonlinearity is critical, so the Cauchy
  problem \eqref{eq.schr2} is globally well posed for small data;
  moreover the blow-up threshold $C_{2,\bt}$ is constant for any $\bt\leq 1$
  and tends to infinity as $\bt\to+\infty$,
\item[iii)] if $n\geq3$ a solution of the Cauchy problem \eqref{eq.schr2}
  exists globally in time provided the initial datum is sufficiently small
  in $L^2\times L^2$.
\end{itemize}

The single equation with a Kerr nonlinearity has been studied also in
bounded domains or on compact manifolds (see, for example, \cite{bgt, bgt2}
and the references therein). The study of coupled nonlinear Schr\"odinger
equations in bounded domains or on compact manifolds should be interesting
in view to extend the results for a single equation.
\eos

The paper is organized in the following way: Section \ref{sec.glob} is devoted
to the proofs of the existence results above, in Section \ref{sec.GN} it is
proved a Gagliardo-Nirenberg inequality (see \eqref{eq.GN}) which is the
fundamental tool to obtain Theorems \ref{thm.existence} and \ref{thm.main}.
Section \ref{sec.BU} deals with a blow-up result which shows the sharpness
of constant $C_{n,\bt}$, while in Section \ref{sec.UT} the proof
of Theorem \ref{thm.main} is completed.

{\bf Acknowledgement\ } We wish to thank P. D'Ancona for the time spent in some 
useful discussions about this work.

\section{Global existence results}\label{sec.glob}

The first part of this work is devoted to the proof of Theorem
\ref{thm.existence}. The theory for the single nonlinear
Schr\"odinger equation \eqref{eq.nls} was developed in \cite{gv1} and \cite{k}; 
the proof of the local well-posedness is a contraction argument based on Strichartz
estimates, and the conservation of both the mass and the energy allows to extend
the local solution globally in time. The fixed point technique also works
in the case of a system, hence problem \eqref{eq.schr1} is
locally well-posed in $H^1$ for $0\leq p\leq2/(n-2)$. We omit here
the straightforward computations, see for example Remarks 4.2.13 and
4.3.4 in \cite{caz}.

Let us now study the conservation laws for system
\eqref{eq.schr1}. Multiplying the equations in \eqref{eq.schr1} by
$\overline\phi$ and $\overline\psi$ respectively, integrating in $x$ and
taking the resulting imaginary parts, we see that
\begin{equation}\label{eq.mass}
  \frac{d}{dt}\|\phi\|_2^2=0,
  \qquad
  \frac{d}{dt}\|\psi\|_2^2=0,
\end{equation}
i.e. the conservation of the masses. Note that these computations make sense
if $\phi,\psi$ are $H^1$ solutions (it is possible to prove
\eqref{eq.mass} also in the case of $L^2$ solutions, following for
example the techniques of \cite{o}).

Now we consider the energy $E(t)$ defined in \eqref{eq.energy}. Let
$(\phi,\psi)$ be a solution to \eqref{eq.schr1}; multiplying the
equations in \eqref{eq.schr1} by $\overline\phi_t$ and
$\overline\psi_t$ respectively, integrating by parts in $x$ and
taking the resulting real parts, we easily obtain the energy
conservation
\begin{equation}\label{eq.energycons}
  E'(t)=0.
\end{equation}
This formal computation needs $H^2$ regularity for $\phi,\psi$, but
\eqref{eq.energycons} makes sense (and can be proved) also for $H^1$ solutions. To prove
this, following exactly the same computations of Ozawa in
\cite{o}, Proposition 2; we omit here the details.

In order to obtain an a priori control on the gradient of the
solutions, we introduce a Gagliardo-Nirenberg inequality (see
Section \ref{sec.GN} below):
\begin{equation}\label{eq.GN}
  \left(\|u\|_{2p+2}^{2p+2}+2\beta\|uv\|_{p+1}^{p+1}
  +\|v\|_{2p+2}^{2p+2}\right)\leq
  C_{n,p,\beta} \left( \|u\|_{2}^{2}+\|v\|_{2}^{2}\right)^{p+1-p\frac n2}
  \left(\|\nabla u\|_{2}^{2}+\|\nabla v\|_{2}^{2}\right)^{p\frac n2},
\end{equation}
that gives the following bound from below:
\beq\label{eq.energybelow}
  E(t)\geq
  \frac{1}{2}(\|\nabla\phi\|_2^2+\|\nabla\psi\|_2^2)
  \left[1-\frac{C_{n,p,\beta}}{p+1}\left(\|\phi\|^{2}_{2}
  +\|\psi\|^2_2 \right)^{p+1-p\frac n2}
  \left(\|\nabla\phi\|^2_2
  +\|\nabla\psi\|^2_2\right)^{p\frac n2-1}\right].
\eeq
If $p<2/n$, we easily see by \eqref{eq.energybelow} that the norms
$\|\nabla\phi\|_2,\|\nabla\psi\|_2$ cannot blow up in a finite time,
because of the conservation of both the mass and the energy; as a consequence,
global well-posedness in $H^1$ is proved in the subcritical range.
The power $p=2/n$ is critical, in the sense that this nonlinearity
is sufficiently high to generate $H^1$ solutions blowing up in a
finite time. On the other hand, also in this case, the smallness
assumption
\begin{equation}\label{eq.criticGN}
  (\|\phi\|^2_2+\|\psi\|^2_2)^{2/n}
  <\frac{p+1}{C_{n,p,\beta}}
\end{equation}
allows by \eqref{eq.energybelow} to obtain the same a priori control
for the gradient in terms of the energy, hence the global existence
in the Theorems \ref{thm.existence} and \ref{thm.main} is proved.
The last part of Theorem \ref{thm.main} is proved
in Section \ref{sec.BU}.

\section{Gagliardo-Nirenberg inequality}\label{sec.GN}

Our next step is to discuss, following the approach of \cite{w},
the behavior of the best constant $C_{n,p,\beta}$ in the Gagliardo-Nirenberg
inequality \eqref{eq.GN}; this will allow us to understand which is the
critical initial level defining the border line between global
well-posedness and blow-up phenomena. This involves the existence of
minimal energy stationary solutions of \eqref{eq.schr1} and allows
us to clarify the concept of {\it ground state}.

Consider the functional
\bdm
  J_{n,p,\beta}(u,v)= \dfrac{\left(\|\nabla u\|_{2}^{2}
  +\|\nabla v\|_{2}^{2}\right)^{pn/2} \left( \|u\|_{2}^{2}+
  \|v\|_{2}^{2}\right)^{p+1-pn/2}}
  {\left(\|u\|_{2p+2}^{2p+2}+2\beta\|uv\|_{p+1}^{p+1}
  +\|v\|_{2p+2}^{2p+2}\right)},\quad u,v\in\huno;
\edm
the infimum of $J_{n,p,\beta}$ on $H^1\times H^1$ is clearly the reciprocal of the best
constant $C_{n,p,\beta}$ in \eqref{eq.GN}. First of all we want
to point out that, for any $u,v \in\huno$ and for any $\mu,\la>0$,
if we set $u_{\mu,\la}(x)=\mu u(\la x)$ and $v_{\mu,\la}(x)=\mu v(\la x)$
it follows
\bdm
\|u_{\mu,\la}\|_{2}^{2}=\mu^{2}\la^{-n}\|u\|_{2}^{2},\quad
\|\nabla u_{\mu,\la}\|_{2}^{2}=\mu^{2}\la^{2-n}\|\nabla u\|_{2}^{2},
\edm
\bdm
\|v_{\mu,\la}\|_{2}^{2}=\mu^{2}\la^{-n}\|v\|_{2}^{2},\quad
\|\nabla v_{\mu,\la}\|_{2}^{2}=\mu^{2}\la^{2-n}\|\nabla v\|_{2}^{2},
\edm
\bdm
\|u_{\mu,\la}\|_{2p+2}^{2p+2}=\mu^{2p+2}\la^{-n}\|u\|_{2p+2}^{2p+2},
\quad
\|v_{\mu,\la}\|_{2p+2}^{2p+2}=\mu^{2p+2}\la^{-n}\|v\|_{2p+2}^{2p+2},
\edm
so that
\bdm
J_{n,p,\beta}(u_{\mu,\la},v_{\mu,\la})=J_{n,p,\beta}(u,v).
\edm

Assume that the infimum of $J_{n,p,\beta}$ is achieved by
$(\tilde{u},\tilde{v})$: since the value of the functional is
invariant with respect to the above scalings, we can assume that the best
constant in \eqref{eq.GN} is achieved by the pair
$(\tilde{u},\tilde{v})$ such that \bdm
\left(\|\tilde{u}\|_{2}^{2}+\|\tilde{v}\|_{2}^{2}\right)=
\left(\|\nabla \tilde{u}\|_{2}^{2}+\|\nabla
\tilde{v}\|_{2}^{2}\right)=1. \edm

Therefore $(\tilde u,\tilde v)$ is a weak solution of the following
system of two weakly coupled elliptic equations
\bdm
\begin{cases}
  -\dfrac{pn}{2}\lapl\tilde{u}+\dfrac{(2-n)p+2}{2}\tilde{u}
  =\dfrac{1}{C_{n,p,\beta}}\left(|\tilde{u}|^{2p}
  +\beta|\tilde{u}|^{p-1}|\tilde{v}|^{p+1}\right)\tilde{u} & \\
  -\dfrac{pn}{2}\lapl \tilde{v}+\dfrac{(2-n)p+2}{2}\tilde{v}
  =\dfrac{1}{C_{n,p,\beta}}\left(|\tilde{v}|^{2p}
  +\beta|\tilde{v}|^{p-1}|\tilde{u}|^{p+1}\right)\tilde{v}, &
\end{cases}
\edm

Now consider the pair $(\tilde{u}_{\mu,\la},\tilde{v}_{\mu,\la})$,
corresponding to the choice of parameters 
\bdm 
\mu=\left(
\dfrac{2}{C_{n,p,\beta}(2p+2-pn)} \right)^{1/2p},\quad \la=\left(
\dfrac{pn}{(2p+2-pn)} \right)^{1/2};
\edm
this pair solves the following elliptic system
\beq\label{eq.ellittico}
\begin{cases}
  -\lapl\tilde{u}_{\mu,\la}+\tilde{u}_{\mu,\la}=\left(|\tilde{u}_{\mu,\la}|^{2p}
  +\beta|\tilde{u}_{\mu,\la}|^{p-1}|\tilde{v}_{\mu,\la}|^{p+1}\right)\tilde{u}_{\mu,\la} & \\
  -\lapl\tilde{v}_{\mu,\la}+\tilde{v}_{\mu,\la}=\left(|\tilde{v}_{\mu,\la}|^{2p}
  +\beta|\tilde{v}_{\mu,\la}|^{p-1}|\tilde{u}_{\mu,\la}|^{p+1}\right)\tilde{v}_{\mu,\la}. &
\end{cases} 
\eeq
Note that the preceding system is variational in nature, so that any (weak) solution
is a critical point of the functional
\bdm
I_{n,p,\beta}(u,v)=\dfrac{1}{2}\left(\|\nabla u\|_{2}^{2}+\|\nabla v\|_{2}^{2}+
\|u\|_{2}^{2}+\|v\|_{2}^{2}\right)-\dfrac{1}{2p+2} \left(\|u\|_{2p+2}^{2p+2}
+\beta\|u v\|_{p+1}^{p+1}+\|v\|_{2p+2}^{2p+2}\right).
\edm

Recently the problem of existence of positive solutions for elliptic systems of this 
kind has been studied by many authors (see, for example, \cite{ac, bw, dfl, lw, mmp, s, y}).
In \cite{mmp} particular attention is given to the existence and some qualitative
properties of the ground state solutions of \eqref{eq.ellittico}: a ground state
solution is a nontrivial solution (i.e. distinct from the pair $(0,0)$) which has the
least critical level. In particular it is possible to prove existence of ground
state solutions for system \eqref{eq.ellittico} solving the following minimization problem
\beq\label{eq.minim}
\inf_{(u,v)\in\mathcal{N}} I_{n,p,\beta}(u,v),
\eeq
where $\mathcal{N}\subset\huno\times\huno$ is the Nehari manifold,
that is
\bdm
  \mathcal{N}=\left\{ (u,v)\neq(0,0):
  \|\nabla u\|_{2}^{2}+\|\nabla v\|_{2}^{2}+\|u\|_{2}^{2}+\|v\|_{2}^{2}
  =\|u\|_{2p+2}^{2p+2} +2\beta\|uv\|_{p+1}^{p+1} \|v\|_{2p+2}^{2p+2}
  \right\} .
\edm
Since $\mathcal{N}$ is a smooth (of class $C^2$) manifold containing all the
nontrivial critical points of the functional, that is all the weak solutions of
\eqref{eq.ellittico}, clearly a ground state solution has to realize the minimum.
We want to point out that $I_{n,p,\beta}$ is bounded from below on $\mathcal{N}$,
so the minimization problem \eqref{eq.minim} is well-posed, moreover it is possible
to prove that any minimizing sequences is compact (up to translations) and that
the minimum is achieved.

Let
\bdm
m_{n,p,\beta}=\inf_{\mathcal{N}} I_{n,p,\beta}
=I_{n,p,\beta}(\tilde{u}_{\mu,\la},\tilde{v}_{\mu,\la})
\edm
be the level of any ground state solution of \eqref{eq.ellittico}; we 
want to prove that there is a direct relation between $C_{n,p,\beta}$ 
and $m_{n,p,\beta}$. Recalling that any critical point of $I_{n,p,\beta}$ 
is a weak solution of \eqref{eq.ellittico}, multiplying \eqref{eq.ellittico} 
by $(\tilde{u}_{\mu,\la},\tilde{v}_{\mu,\la})$ and integrating on $\rn$ 
we obtain that
\beq\label{eq.formadebole}
\begin{array}{l}
  \|\nabla\tilde{u}_{\mu,\la}\|_{2}^{2}+\|\tilde{u}_{\mu,\la}\|_{2}^{2}
  =\|\tilde{u}_{\mu,\la}\|_{2p+2}^{2p+2}
  +\beta\|\tilde{u}_{\mu,\la}\tilde{v}_{\mu,\la}\|_{p+1}^{p+1},  \\ \\
  \|\nabla\tilde{v}_{\mu,\la}\|_{2}^{2}+\|\tilde{v}_{\mu,\la}\|_{2}^{2}
  =\|\tilde{v}_{\mu,\la}\|_{2p+2}^{2p+2}
  +\beta\|\tilde{u}_{\mu,\la}\tilde{v}_{\mu,\la}\|_{p+1}^{p+1}.
\end{array}
\eeq
Moreover in this case, Pohozaev identity reads
\beq\label{eq.pohozaev}
  \begin{array}{l}
  \dfrac{n-2}{2}\left(\|\nabla \tilde{u}_{\mu,\la}\|_2^2
  +\|\nabla \tilde{v}_{\mu,\la}\|_2^2\right)
  +\dfrac{n}{2} \left(\|\tilde{u}_{\mu,\la}\|_2^{2}
  +\|\tilde{v}_{\mu,\la}\|_2^2\right)  \\
  \qquad\qquad=\dfrac{n}{2p+2}\left(\|\tilde{u}_{\mu,\la}\|_{2p+2}^{2p+2}
  +2\beta\|\tilde{u}_{\mu,\la}\tilde{v}_{\mu,\la}\|_{p+1}^{p+1}
  +\|\tilde{v}_{\mu,\la}\|_{2p+2}^{2p+2}\right).
\end{array}
\eeq

Putting together the above identities we have that
\bdm
  \mu^{2}\la^{2-n}\left(\|\nabla\tilde{u}\|_{2}^{2}+\|\nabla\tilde{v}\|_{2}^{2}
  \right)= \left(\|\nabla \tilde{u}_{\mu,\la}\|_{2}^{2}
  +\|\nabla \tilde{v}_{\mu,\la}\|_{2}^{2}\right)= n m_{n,p,\beta},
\edm
\bdm
  \mu^{2}\la^{-n}\left(\|\tilde{u}\|_{2}^{2}+\|\tilde{v}\|_{2}^{2}\right)
  =\left(\|\tilde{u}_{\mu,\la}\|_{2}^{2}+\|\tilde{v}_{\mu,\la}\|_{2}^{2}\right)
  =\left(2-n+\dfrac{2}{p}\right) m_{n,p,\beta},
\edm
\bdm
\begin{array}{l}
  \mu^{2p+2}\la^{-n}\left(\|\tilde{u}\|_{2p+2}^{2p+2}
  +2\beta\|\tilde{u}\tilde{v}\|_{p+1}^{p+1}
  +\|\tilde{v}\|_{2p+2}^{2p+2}\right) \\
  \qquad=\left(\|\tilde{u}_{\mu,\la}\|_{2p+2}^{2p+2}
  +2\beta\|\tilde{u}_{\mu,\la}\tilde{v}_{\mu,\la}\|_{p+1}^{p+1}
  +\|\tilde{v}_{\mu,\la}\|_{2p+2}^{2p+2}\right)
  =\dfrac{2p+2}{p}m_{n,p,\beta}.
\end{array}
\edm

All the above calculations imply that the following equalities hold
\beq\label{eq.mec}
\dfrac{1}{C_{n,p,\beta}}=J_{n,p,\beta}(\tilde{u},\tilde{v})
= J_{n,p,\beta}(\tilde{u}_{\mu,\la},\tilde{v}_{\mu,\la})
=m_{n,p,\beta}^p\dfrac{n^{pn/2}(2p+2-pn)^{p+1-pn/2}}{2(p+1)p^{p-pn/2}}.
\eeq
Note that, in the critical case $p=2/n$, \eqref{eq.mec} becomes
\beq\label{eq.meccr}
\dfrac{1}{C_{n,2/n,\beta}}=2^{2/n}\dfrac{n}{n+2}m_{n,2/n,\beta}^{2/n}
=\dfrac{n}{n+2}\left(\|\tilde{u}_{\mu,\la}\|_{2}^{2}
+\|\tilde{v}_{\mu,\la}\|_{2}^{2}\right)^{2/n}.
\eeq

The arguments above, in particular \eqref{eq.mec}, show that a
minimum point of $J_{n,p,\beta}$, through a suitable scaling, has to
correspond to a ground state solution of \eqref{eq.ellittico} (or to a
least energy nontrivial critical point of $I_{n,p,\beta}$). Now,
since in \cite{mmp} it is proved the existence of ground state
solutions to \eqref{eq.ellittico}, we have obtained the existence of
a minimum point for the functional $J_{n,p,\beta}$; this shows that
inequality \eqref{eq.GN} is sharp and that there exists at least a
pair of functions for which equality holds. More generally we have 
proved that the functionals $J_{n,p,\beta}$ and $I_{n,p,\beta}$ possess
the same number of critical values.

The validity of inequality \eqref{eq.GN} follows by the above arguments.

\section{Blow-up results}\label{sec.BU}

In view to prove the sharpness of the constant $C$ inthe statement
of Theorem \ref{thm.main}, we introduce (following \cite{gl} and \cite{rn})
another physically relevant quantity, that plays a crucial
role in the analysis of blow-up phenomena: the {\it variance}
$V(t)$, which is defined by
\begin{equation}\label{eq.variance}
  V(t)=\int |x|^2|\phi(t,x)|^2\,dx
  +\int |x|^2|\psi(t,x)|^2\,dx.
\end{equation}
As in the case of a single Schr\"odinger equation, we will prove
a relation between the time behavior of $V$ and that of the $H^1$-norm of
the solutions: as we will see in the following, the precise calculation
of the first and second derivatives of $V$ in terms of the solutions
of \eqref{eq.schr1} is the main tool for the description of the
blow-up (see for example \cite{caz} for a proof in the case of a
single equation).

More precisely, we prove the following Lemma:
\begin{lemma}\label{lem.variance}
  Let $(\phi,\psi)$ be a solution of system \eqref{eq.schr1} on an
  interval $I=(-t_1,t_1)$; then, for each $t\in I$, the variance
  satisfies the following
  identities:
  \begin{align}\label{eq.Vprimo}
    V'(t)= & 4\Im\int\left[\bigl(x\cdot
    \nabla\phi\bigr)
    \overline{\phi}
    +\bigl(x\cdot\nabla\psi\bigr)
    \overline{\psi}\right]\,dx,
    \\
    V''(t)= &
    8\int\bigl(|\nabla\phi|^2+|\nabla\psi|^2\bigr)\,dx
    -\frac{4np}{p+1}\int
    \bigl(|\phi|^{2p+2}+2\beta |\phi\psi|^{p+1}+|\psi|^{2p+2}\bigr)\,dx.
    \label{eq.Vsecondo}
  \end{align}
\end{lemma}

\dimo We introduce the following notations:
\begin{align*}
  z= & (z^1,\dots,z^n)\in\C^n;
  \\
  z\cdot w= & \sum_iz^iw^i,\qquad z,w\in\C^n;
  \\
  u_i= & \frac{\partial u}{\partial x_i},
  \qquad u:\R^n\to\C.
\end{align*}

Multiplying the equations in \eqref{eq.schr1} by $2\overline\phi$
and $2\overline\psi$ respectively, and taking the resulting
imaginary parts, we obtain
\begin{align}
  \frac{\partial}{\partial t}|\phi|^2= &
  -2\Im(\overline\phi\Delta\phi)=
  -2\nabla\cdot(\Im\overline\phi\nabla\phi),
  \label{eq.1}
  \\
  \frac{\partial}{\partial t}|\psi|^2= &
  -2\Im(\overline\psi\Delta\psi)=
  -2\nabla\cdot(\Im\overline\psi\nabla\psi).
  \label{eq.2}
\end{align}
Now, multiplying \eqref{eq.1} and \eqref{eq.2} by $|x|^2$,
and integrating by parts in $x$, we immediately obtain
\eqref{eq.Vprimo}.

In order to prove \eqref{eq.Vsecondo}, let us multiply the
equations in \eqref{eq.schr1} by $2(x\cdot\nabla\overline\phi)$
and $2(x\cdot\nabla\overline\psi)$ respectively, let us integrate
in $x$ and sum the equations for the real parts, to get:
\begin{align*}
  0= &
  2\Re\int i\bigl[(x\cdot\nabla\overline\phi)\phi_t
  (x\cdot\nabla\overline\psi)\psi_t\bigr]\,dx
  +2\Re\int\bigl[(x\cdot\nabla\overline\phi)
  \Delta\phi+(x\cdot\nabla\overline\psi)
  \Delta\psi\bigr]\,dx
  \\
  & +2\Re\int\bigl[(x\cdot\nabla\overline\phi)
  (|\phi|^{2p}+\beta|\psi|^{p+1}
  |\phi|^{p-1})\phi+\bigr.
  \\
  & \bigl.+
  (x\cdot\nabla\overline\psi)
  (|\psi|^{2p}+\beta|\phi|^{p+1}
  |\psi|^{p-1})\psi\bigr]\,dx.
\end{align*}
We rewrite the last identity in the form
\begin{equation}\label{eq.eq}
  \hbox{I}=\hbox{II}+\hbox{III},
\end{equation}
where
\begin{align*}
  \hbox{I}= & 2\Re\int i\bigl[(x\cdot\nabla\overline\phi)\phi_t
  (x\cdot\nabla\overline\psi)\psi_t\bigr]\,dx, 
  \\
  \hbox{II}= & -2\Re\int\bigl[(x\cdot\nabla\overline\phi)
  \Delta\phi+(x\cdot\nabla\overline\psi)
  \Delta\psi\bigr]\,dx, 
  \\
  \hbox{III}= & -2\Re\int\bigl[(x\cdot\nabla\overline\phi)
  (|\phi|^{2p}+\beta|\psi|^{p+1}
  |\phi|^{p-1})\phi\bigr. 
  \\
  & \bigl.+(x\cdot\nabla\overline\psi)
  (|\psi|^{2p}+\beta|\phi|^{p+1}
  |\psi|^{p-1})\psi\bigr]\,dx.
\end{align*}
For the first term, we have
\begin{equation*}
  \hbox{I}=
  -\Re\int i\sum_j\left(x^j\overline\phi_j\phi_t
  -x^j\phi_j\overline\phi_t+
  x^j\overline\psi_j\psi_t
  -x^j\psi_j\overline\psi_t\right)\,dx,
\end{equation*}
which can be written in the form
\begin{align*}
  \hbox{I}= &
  \Re\int i\sum_jx^j\left[(\overline\phi_j\phi)_t-
  (\phi\overline\phi_t)_j
  +(\overline\psi_j\psi)_t-
  (\psi\overline\psi_t)_j\right]\,dx
  \\
  = &
  \frac{d}{dt}\Re\int i\left[
  (x\cdot\nabla\overline\phi)
  \phi+(x\cdot\nabla\overline\psi)
  \psi\right]\,dx
  +n\Re\int i(\phi\overline\phi_t+
  \psi\overline\psi_t)\,dx.
\end{align*}
Now we evaluate the last equality using the
equations in \eqref{eq.schr1}, obtaining
\begin{align}\label{eq.Afin}
  \hbox{I} &
  =\frac{d}{dt}\Im\int\left[
  (x\cdot\nabla\phi)
  \overline\phi+(x\cdot\nabla\psi)
  \overline\psi\right]\,dx
  -n\int\left(|\nabla\phi|^2+
  |\nabla\psi|^2\right)\,dx
  \\
  &
  \qquad+n\int
  \left[(|\phi|^{2p}+\beta|\psi|^{p+1}
  |\phi|^{p-1})|\phi|^2+
  (|\psi|^{2p}+\beta|\phi|^{p+1}
  |\psi|^{p-1})|\psi|^2\right]\,dx
  \nonumber \\
  &
  =\frac{d}{dt}\Im\int\left[
  (x\cdot\nabla\phi)
  \overline\phi+(x\cdot\nabla\psi)
  \overline\psi\right]\,dx
  -n\int\left(|\nabla\phi|^2+
  |\nabla\psi|^2\right)\,dx
  \nonumber \\
  &
  \qquad+n\int
  \left(|\phi|^{2p+2}+2\beta|\psi\phi|^{p+1}
  +|\psi|^{2p+2}\right)\,dx.
  \nonumber
\end{align}
A multiple integration by parts in $\hbox{II}$ gives the
Pohozaev identity
\begin{equation}\label{eq.Bfin}
  \hbox{II}=(2-n)\int\left(|\nabla\phi|^2+
  |\nabla\psi|^2\right)\,dx.
\end{equation}
As for the term $\hbox{III}$, we write it by components:
\begin{align}\label{eq.C2}
  \hbox{III}= & -\sum_j\int\left\{x^j\left[|\phi|^{2p}
  (2\Re\overline\phi_j\phi)
  +|\psi|^{2p}
  (2\Re\overline\psi_j\psi)\right]\right.
  \\
  & \left.+
  \beta x^j\left[|\phi|^{p-1}|\psi|^{p+1}
  (2\Re\overline\phi_j\phi)+
  |\psi|^{p-1}|\phi|^{p+1}
  (2\Re\overline\psi_j\psi)\right]\right\}\,dx.
  \nonumber
\end{align}
Observe that
\begin{align*}
  & |\phi|^{2p}
  (2\Re\overline\phi_j\phi)
  +|\psi|^{2p}
  (2\Re\overline\psi_j\psi)=
  \frac{1}{p+1}
  \left(|\phi|^{2p+2}_j+|\psi|^{2p+2}_j\right),
  \\
  & |\phi|^{p-1}|\psi|^{p+1}
  (2\Re\overline\phi_j\phi)+
  |\psi|^{p-1}|\phi|^{p+1}
  (2\Re\overline\psi_j\psi)=
  \frac{2\beta}{p+1}\left(
  |\phi|^{p+1}|\psi|^{p+1}\right)_j;
\end{align*}
hence, integrating by parts in \eqref{eq.C2} we have
\begin{equation}\label{eq.Cfin}
  \hbox{III}=\frac{n}{p+1}
  \int\left(|\phi|^{2p+2}+|\psi|^{2p+2}
  +2\beta |\phi|^{p+1}|\psi|^{p+1} \right)\,dx.
\end{equation}
Finally, recollecting \eqref{eq.eq}, \eqref{eq.Afin},
\eqref{eq.Bfin}, \eqref{eq.Cfin} and \eqref{eq.Vprimo}, we complete
the proof of \eqref{eq.Vsecondo}. \finedim

\bos
Note that \eqref{eq.Vsecondo} can be rewritten, recalling the definition 
of $E$, in the following equivalent form
\beq\label{eq.varianza}
  V''(t)= 16 E(t) -8\frac{np-2}{2p+2}\int
  \bigl(|\phi|^{2p+2}+2\beta |\phi\psi|^{p+1}+|\psi|^{2p+2}\bigr)\,dx.
\eeq

In the critical case $p=2/n$ the equation above reduces to
\bdm
  V''(t)= 16 E(t);
\edm
hence the variance $V$ of any solution of \eqref{eq.schr1} with negative 
initial energy vanish in a finite time. For each $h:\R^n\to\C$ we can estimate
\begin{equation*}
  \|h\|_{L^2}^2=\||h|^2\|_{L^1}\leq\||x|h\|_{L^2}\left\|\frac{h}{|x|}\right\|_{L^2},
\end{equation*}
by Cauchy-Schwartz inequality. As a consequence of the standard Hardy's inequality 
we obtain
\begin{equation*}
  \|h\|_{L^2}^2\leq\||x|h\|_{L^2}\|\nabla h\|_{L^2}.
\end{equation*}
Applying the last inequality to any solution of of \eqref{eq.schr1} with negative 
initial energy, since the mass is conserved and the variance vanish in a finite 
time the $L^2$ norm of the gradient needs necessarely to blow up in a finite time.
\eos

\bos
Consider the following pair
\bdm
\frac{e^{-i\tfrac{|x|^2-4}{4(1-t)}}}{(1-t)^{n/2}}
\left(U\left(\frac{x}{1-t}\right),V\left(\frac{x}{1-t}\right)\right),
\edm
where $(U,V)$ is a ground state solution of \eqref{eq.ellittico}.
This is an explicit example of a blow-up solution, which shows that
the constant $C_{n,\bt}$ in Theorem \ref{thm.main} is sharp. 
Indeed the above pair solves \eqref{eq.schr1}
with initial value $e^{-i(|x|^2-4)/4}(U(x),V(x))$
which attains the critical blow-up threshold.
\eos

\section{On the blow-up threshold}\label{sec.UT}

If $p=2/n$, we have obtained the following characterization of the
blow-up threshold
\beq\label{eq.charact}
\dfrac{1}{C_{n,2/n,\beta}}=\dfrac{1}{C_{n,\beta}}=\dfrac{n}{n+2}
\left(\|U\|_{2}^{2}+\|V\|_{2}^{2}\right)^{2/n},
\eeq
where $(U,V)$ is a ground state solution of \eqref{eq.ellittico}.
In order to prove Theorem \ref{thm.main} we have to estimate the
quantities involved in \eqref{eq.charact}.

In \cite{mmp} (see Theorem $2.5$) it is proved that, if $\bt< 2^{2/n}-1$, then 
any ground state of the elliptic system \eqref{eq.ellittico} is a scalar
function, that is one of the components of the ground state solution is zero.
So we can assume, without loss of generality, that the ground state is
$(z,0)$, where $z\in H^1$ is the unique ground state solution (see \cite{w}) of the
equation
\beq\label{eq.single}
-\lapl z+ z=|z|^\frac{4}{n} z.
\eeq
This implies that the constant $C_{n,\beta}=C_{n}$ depends only on $n$ for any
$\bt\leq 2^{2/n}-1$, since the coupling parameter $\bt$ now does not play a role
in the problem of selecting the ground state solution.
Moreover $C_{n}$ is exactly the blow-up threshold for a single nonlinear
Schr\"odinger equation, introduced and numerically computed in \cite{w}.

If $\bt\geq 2^{2/n}-1$, $C_{n,\beta}$ depends on $n$ and $\bt$ and
its expression is unknown, but we can estimate it using a suitable
test-pair. Let $\hat{z}$ be the unique positive ground state
solution of 
\bdm 
-\lapl \hat{z}+ \hat{z}=(1+\bt)|\hat{z}|^\frac{4}{n}\hat{z}, 
\edm 
it is easy to see that the pair $(\hat{z},\hat{z})$ is a
positive solution of \eqref{eq.ellittico} for any $\bt$;
and the following inequality holds 
\bdm
\dfrac{1}{C_{n,\beta}}=\dfrac{n}{n+2}\left(\|U\|_{2}^{2}+\|V\|_{2}^{2}\right)^{2/n}
\leq \dfrac{n}{n+2}\left(2\|\hat{z}\|_{2}^{2}\right)^{2/n}. 
\edm
Clearly we have an inequality since $(\hat{z},\hat{z})$ could not be
a ground state solution of \eqref{eq.ellittico}.

Using the scaling as in Section \ref{sec.GN} we can estimate $C_{n,\beta}$
with $C_{n}$. Recalling that $z$ is the ground state solution of
\eqref{eq.single} and noticing that the $L^2$ norm of $z$ is related
to $C_n$ (see $(1.3)$ in \cite{w} and also \eqref{eq.meccr}) we obtain that
\bdm
\|\hat{z}\|_{4/n}^{2}=\dfrac{\|z\|_{4/n}^{2}}{1+\bt}
=\dfrac{n+2}{n(1+\bt)C_n}.
\edm
Collecting the inequalities above we have
\bdm
C_{n,\beta}\geq \left(\dfrac{1+\bt}{2^{2/n}}\right) C_{n},
\edm
so that the claim is proved.\finedim

Now we give some concluding remarks.

\bos
Consider, on the space $H^1\times H^1$, the functionals
\bdm
  J^{\om}_{n,p,\beta}(u,v)= \dfrac{\left(\|\nabla u\|_{2}^{2}
  +\|\nabla v\|_{2}^{2}\right)^{pn/2} \left( \om_{1}\|u\|_{2}^{2}
  +\om_{2}\|v\|_{2}^{2}\right)^{p+1-pn/2}}
  {\left(\|u\|_{2p+2}^{2p+2}+2\beta\|uv\|_{p+1}^{p+1}
  +\|v\|_{2p+2}^{2p+2}\right)},\quad u,v\in\huno;
\edm
with $\om_{1},\om_{2}$ positive real numbers and $p,\bt$ as above
(this corresponds to looking for ground state solutions in the more 
general form $(\phi,\psi)=\left(e^{i\om_1 t}u(x),e^{i\om_2 t}v(x)\right)$).

The infimum of $J^{\om}_{n,p,\beta}$ on $H^1\times H^1$ is the reciprocal
of the best constant $C^{\om}_{n,p,\beta}$ in the following inequality
\beq\label{eq.GNomega}
  \left(\|u\|_{2p+2}^{2p+2}+2\beta\|uv\|_{p+1}^{p+1}
  +\|v\|_{2p+2}^{2p+2}\right)\leq
  C^{\om}_{n,p,\beta}\left(\om_{1}\|u\|_{2}^{2}
  +\om_{2}\|v\|_{2}^{2}\right)^{p+1-p\frac n2}
  \left(\|\nabla u\|_{2}^{2}+\|\nabla v\|_{2}^{2}\right)^{p\frac n2},
\eeq
which is a little generalization of \eqref{eq.GN}. Arguing as in Section 
\ref{sec.GN} we can prove that the constant $C^{\om}_{n,p,\bt}$ is achived, 
up to scaling, by a ground state solution of the following elliptic system
\beq\label{eq.ellitticomega}
\begin{cases}
  -\lapl u+\om_{1}u=\left(|u|^{2p}+\beta|u|^{p-1}|v|^{p+1}\right)u & \\
  -\lapl v+\om_{2}v=\left(|v|^{2p}+\beta|v|^{p-1}|u|^{p+1}\right)v. &
\end{cases}
\eeq
The existence and some qualitative properties of the ground state solutions
of \eqref{eq.ellitticomega} are discussed in \cite{mmp}.
Inequalities \eqref{eq.GNomega} can be used in order to prove the above
existence results, Theorems \ref{thm.existence} and \ref{thm.main}.

It is an interesting open question if these inequalities can improve
the existence results stated in Section \ref{sec.introd}.Note that
it is possible to exhibit a blow-up solution of \eqref{eq.schr1}
starting from any ground state solution of \eqref{eq.ellitticomega},
for any choice of two positive real numbers $\om_{1},\om_{2}$,
we can consider the following pair
\beq\label{eq.buomega}
(1-t)^{-n/2} \left(e^{-i\tfrac{|x|^2-4\om_{1}}{4(1-t)}}U_{\om}\left(\frac{x}{1-t}\right),
e^{-i\tfrac{|x|^2-4\om_{2}}{4(1-t)}}V_{\om}\left(\frac{x}{1-t}\right)\right),
\eeq
where $(U_{\om},V_{\om})$ is a ground state solution of \eqref{eq.ellitticomega}.
\eos

\bos
Note that it is possible to extend this argument to systems with
more than two nonlinear Schr\"odinger equations, using some results about
the elliptic counterpart contained in \cite{ac}, Section $6$.
\eos


\end{document}